\documentclass{amsart}
\usepackage{graphicx}

\DeclareMathOperator{\Cos}{Cos}

\newcommand{\C}{{\mathbb C}}

\begin{document}

\title{A Gronwall-type Trigonometric Inequality}
\author{A.~G.~Smirnov}
\address{I.~E.~Tamm Theory Department, P.~N.~Lebedev
Physical Institute, Leninsky prospect 53, Moscow 119991, Russia}
\email{smirnov@lpi.ru}

\begin{abstract}
We prove that the absolute value of the $n$th derivative of $\cos(\sqrt{x})$ does not exceed $n!/(2n)!$ for all $x>0$
and $n = 0,1,\ldots$ and obtain a natural generalization of this inequality involving the analytic continuation of $\cos(\sqrt{x})$.
\end{abstract}

\maketitle

In 1913, Gronwall~\cite{Gronwall} proposed to prove that
\begin{equation}\label{0}
\left|\frac{d^n}{dx^n}\frac{\sin x}{x}\right| \leq \frac{1}{n+1}
\end{equation}
for all real $x$ and all $n = 0,1,\ldots$. A short and elegant proof of this inequality was given by Dunkel~\cite{Dunkel} in 1920. (Another proof in~\cite{Dunkel} provided by Uhler was extremely cumbersome.)
In this note, another trigonometric inequality of a similar type is established: we shall show that
\begin{equation}\label{1}
\left|\frac{d^n}{dx^n}\cos(\sqrt{x})\right| \leq \frac{n!}{(2n)!}, \quad x > 0,
\end{equation}
for all $n = 0,1,\ldots$. While this statement can hardly be directly derived from~(\ref{0}), its demonstration given below bears a strong resemblance to Dunkel's proof of~(\ref{0}).

Our motivation for estimating the derivatives of $\cos(\sqrt{x})$ comes from the study of one-dimensional Schr\"odinger operators of the form
\begin{equation}\label{inv}
-\partial^2_r + \frac{\alpha}{r^2},\quad r>0,
\end{equation}
where $\alpha$ is a real parameter. This system (which is the simplest one-particle variant of the well-known Calogero model~\cite{Calogero}) exhibits a phase transition~\cite{GTV} at $\alpha = -1/4$: the self-adjoint realizations of differential expression~(\ref{inv}) have infinitely many eigenvalues for $\alpha< - 1/4$ and at most one eigenvalue for $\alpha\geq -1/4$. In~\cite{Smirnov}, we used the method of singular Titchmarsh--Weyl $m$-functions~\cite{KST} to construct eigenfunction expansions for~(\ref{inv}) in a neighborhood of the critical point $\alpha = -1/4$ and proved that the corresponding spectral measures depend continuously on $\alpha$. Inequality~(\ref{1}) (or, more precisely, its generalization~(\ref{13}) below) provides us with a tool for estimating the $\alpha$-derivatives of the $m$-functions for this model. This makes it possible to prove that the spectral measures are actually smooth in $\alpha$ near the critical point. (A detailed treatment of this problem will be given in a forthcoming paper.)

\begin{figure}
  \includegraphics[width=\linewidth]{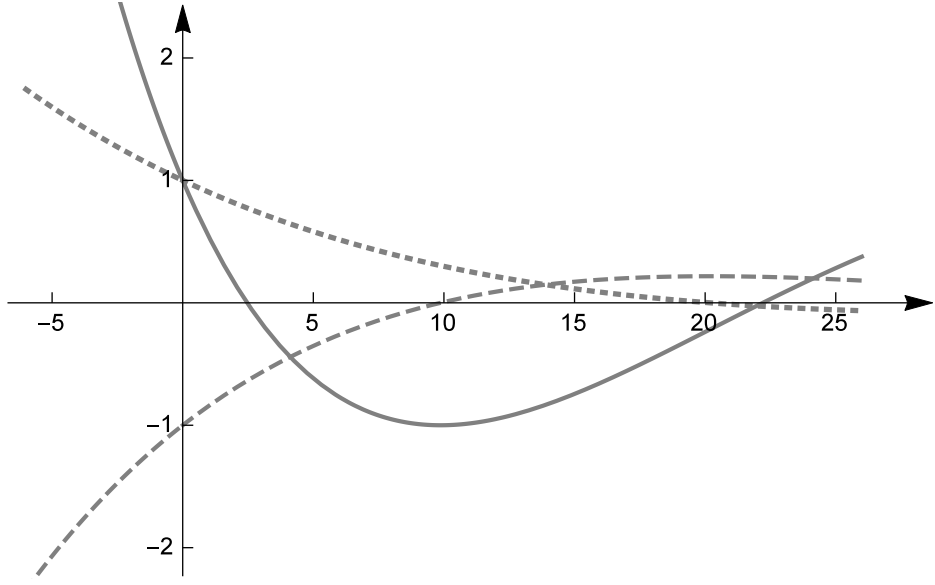}
  \caption{Solid, dashed, and dotted lines correspond to the functions
  $n!^{-1}(2n)!\Cos^{(n)}x$
  for $n = 0$, $1$, and $2$  respectively. According to~(\ref{5}), these functions are equal to $\pm 1$ at $x = 0$.}
  \label{f1}
\end{figure}

Let the entire analytic function $\Cos$ be defined by the relation
\begin{equation}\label{2}
\Cos z = \sum_{k = 0}^\infty \frac{(-1)^k z^k}{(2k)!},\quad z\in\C.
\end{equation}
Then $\Cos(z^2) = \cos z$ for all $z\in\C$. Hence, $\Cos z$ is the analytic continuation of $\cos(\sqrt{x})$ from the positive half-axis to the entire complex plane, and we have
\begin{equation}\label{2a}
\Cos x = \Cos(-|x|) = \cos\left(i \sqrt{|x|}\right) = \cosh\left(\sqrt{|x|}\right),\quad x\leq 0.
\end{equation}
(The graphs of $\Cos x$ and its first two derivatives are shown in Figure~\ref{f1}.)

It follows immediately from~(\ref{2}) that
\begin{equation}\label{5}
\Cos^{(n)}(0) = (-1)^n \frac{n!}{(2n)!},
\end{equation}
where $\Cos^{(n)}z$ denotes the $n$th derivative of $\Cos z$.
This equality implies that the constant in the right-hand side of~(\ref{1}) is optimal. For all $z\in\C$, we have
\begin{equation}\label{7}
(\Cos z)^2 + 4 z (\Cos' z)^2 = 1.
\end{equation}
Indeed, since
\begin{equation}\label{8}
\Cos' x = -\frac{\sin(\sqrt{x})}{2\sqrt{x}},\quad x>0,
\end{equation}
equation~(\ref{7}) holds on the positive half-axis. By the uniqueness theorem for analytic functions, (\ref{7}) remains valid for all $z\in \C$. Differentiating~(\ref{7}), we obtain
\begin{equation}\nonumber
\Cos z + 2 \Cos' z + 4 z\Cos'' z = 0,\quad z\in\C.
\end{equation}
By induction on $n$, this implies that
\begin{equation}\label{10}
\Cos^{(n-1)}z + (4n - 2)\Cos^{(n)}z + 4z \Cos^{(n+1)}z = 0,\quad z\in\C,
\end{equation}
for all $n \geq 1$. It follows from~(\ref{8}) that both $\Cos' x$ and $\Cos''x$ tend to zero as $x\to\infty$. Using induction on $n$, we deduce from~(\ref{10}) that
\begin{equation}\label{11}
\lim_{x\to\infty}\Cos^{(n)}x = 0
\end{equation}
for all $n \geq 1$.

We now prove~(\ref{1}) by induction on~$n$. Clearly, (\ref{1}) is valid for $n = 0$. Let $n \geq 1$ and suppose that (\ref{1}) holds for all derivatives of orders less than $n$. By~(\ref{11}), the maximum value of $|\Cos^{(n)}x|$ on $[0,\infty)$ is attained at some $x_0\geq 0$. If $x_0 = 0$, then the required estimate is ensured by~(\ref{5}). If $x_0>0$, then we have $\Cos^{(n+1)}x_0 = 0$. By (\ref{10}) and the induction hypothesis, it follows that
\[
|\Cos^{(n)}x_0| = \frac{|\Cos^{(n-1)}x_0|}{4n - 2} \leq \frac{(n-1)!}{(2n - 2)!(4n-2)} = \frac{n!}{(2n)!}.
\]
This completes the proof of~(\ref{1}).

It should be noted that inequality~(\ref{1}) is strict for $n\geq 1$, i.e., $|\Cos^{(n)}x| < n!/(2n)!$ for all $x>0$ and $n\geq 1$. Indeed, by~(\ref{8}), this is true for $n = 1$. Suppose $n>1$ and $|\Cos^{(n)}x_0| = n!/(2n)!$ for some $x_0>0$. Then $\Cos^{(n+1)} x_0 = 0$ and it follows from~(\ref{10}) that $|\Cos^{(n-1)}x_0| = (n-1)!/(2n-2)!$. Iterating this argument, we find that $|\Cos^{(1)}x_0| = 1/2$ in contradiction to~(\ref{8}), and our claim is proved.

Using~(\ref{1}) and the power series expansion for $\Cos z$, it is easy to estimate $\Cos^{(n)}x$ on every interval $[a,\infty)$ with $a\leq 0$. More precisely, we shall establish the inequality
\begin{equation}\label{12}
|\Cos^{(n)} x| \leq \frac{n!(2m)!}{(2n)!m!} |\Cos^{(m)}a|,\quad x\geq a,
\end{equation}
for every $a\leq 0$ and all nonnegative integer numbers $n$ and $m$ such that $m\leq n$. For $a = 0$, this estimate follows immediately from~(\ref{1}) and~(\ref{5}). If $a = 0$ and $m = 0$, then (\ref{12}) coincides with~(\ref{1}). In view of~(\ref{2a}), substituting $m = 0$ in~(\ref{12}) yields
\begin{equation}\label{13}
|\Cos^{(n)} x| \leq \frac{n!}{(2n)!} \cosh\left(\sqrt{|a|}\right),\quad x\geq a.
\end{equation}
We now prove~(\ref{12}). By differentiating~(\ref{2}), we find that
\begin{equation}\nonumber
\Cos^{(n)} z = (-1)^n\sum_{k = 0}^\infty c(n,k)(-z)^k,\quad c(n,k) = \frac{(k+n)!}{k!(2k+2n)!},
\end{equation}
for all $z\in \C$. This implies that
\begin{equation}\label{14}
|\Cos^{(n)}x| = \sum_{k = 0}^\infty c(n,k)|x|^k,\quad x\leq 0.
\end{equation}
It follows that $|\Cos^{(n)}x|$ is decreasing on the negative half-axis and, therefore,
\begin{equation}\label{15}
|\Cos^{(n)}x| \leq |\Cos^{(n)}a|
\end{equation}
for every $x\in [a,0]$. For $x>0$, we have $|\Cos^{(n)}x| \leq |\Cos^{(n)}(0)|$ by~(\ref{1}) and~(\ref{5}) and, hence, (\ref{15}) holds for all $x\geq a$. It easily follows by induction on $n$ that
\begin{equation}\nonumber
c(n,k) \leq \frac{n!(2m)!}{(2n)!m!} c(m,k)
\end{equation}
for all $k=0,1,\ldots$. In view of~(\ref{14}), this implies that (\ref{12}) holds for $x = a$. Combining (\ref{12}) for $x = a$ with~(\ref{15}) yields inequality~(\ref{12}) for arbitrary $x\geq a$.

In conclusion, we return briefly to inequality~(\ref{0}). In fact, (\ref{0}) is an easy consequence of the obvious formula
\begin{equation}\label{repr}
\frac{\sin x}{x} = \int_0^1 \cos\left(tx\right)\,dt.
\end{equation}
Indeed, differentiating~(\ref{repr}) yields
\[
\frac{d^n}{dx^n}\frac{\sin x}{x} = \int_0^1 t^n \cos\left(tx + \frac{n}{2}\pi\right)\,dt,
\]
whence (\ref{0}) follows immediately. Performing the change of variables $t\to y/x$, we obtain the identity
\[
\frac{d^n}{dx^n}\frac{\sin x}{x} = \frac{1}{x^{n+1}}\int_0^x y^n \sin\left(y + \frac{n+1}{2}\pi\right)\,dy,
\]
which was proposed by Gronwall~\cite{Gronwall1} and was mentioned to imply~(\ref{0}) in Section~3.4.24 of~\cite{Mitrinovic}. (It is interesting to note that Gronwall himself failed to link this identity to~(\ref{0}).) It is an open question whether inequality~(\ref{1}) can be derived from some similar integral representations for $\cos(\sqrt{x})$ and its derivatives.

\section*{Acknowledgment} The author is grateful to the referee for pointing out that inequality~(\ref{0}) easily follows from formula~(\ref{repr}).

\end{document}